\documentclass[12pt]{article}
\usepackage{CJK}

\usepackage{xcolor}         
\usepackage{latexsym}
\usepackage{amsmath}
\usepackage{amssymb}
\usepackage{times}
\usepackage{graphicx}
\usepackage{epsfig}
\usepackage{array}
\usepackage{flafter}
\usepackage[page]{appendix}
\usepackage{listings}
\usepackage{algorithm,algpseudocode}
\usepackage{ulem}
\usepackage{verbatim}
\usepackage{hyperref}
\usepackage{qtree}
\usepackage{mathrsfs,amsmath}

\graphicspath{ {figures/} }

\newtheorem{definition}{Definition}
\newtheorem{theorem}{Theorem}
\newtheorem{lemma}{Lemma}
\newtheorem{ExampleDef}{Example}

\newcommand{\Example}[3]{
  \begin{list}{}{
      \setlength{\leftmargin}{1em}} 
    \item                           
    \small                          
    \begin{ExampleDef} \rm          
      {\bf \hspace{-1ex}: #1}       
      #2                                
      \hfill {\large \boldmath $\Box$}  
      \label{ex:#3}                      
    \end{ExampleDef}
  \end{list}}

\title{
Natural Metrics in Contraction Analysis}
\date{}

\begin{document}

\maketitle
\vspace*{-2.5cm}

\begin{center}
{\large Winfried Lohmiller and Jean-Jacques Slotine \par}
{Nonlinear Systems Laboratory \\
Massachusetts Institute of Technology \\
Cambridge, Massachusetts, 02139, USA\\
{\sl \{wslohmil, jjs\}@mit.edu} \par}
\end{center}
\vspace{0.1cm}
\begin{abstract}
Contraction analysis establishes exponential incremental convergence of a nonlinear system by solving a linear matrix inequality for a contraction metric, and has become a standard resource for solving problems in nonlinear control and estimation. This paper shows that, for a general nonlinear system, a contraction metric can be systematically derived by rewriting the system dynamics as a complex natural gradient dynamics. In this form, the variational dynamics can be modally decomposed with quadratic geodesic coordinates, and exact exponential convergence rates can be computed analytically. 

Specializing the results above to Hamiltonian systems shows that differential lengths of general Hamiltonian dynamics correspond to exact complex analytic exponential functions, whose eigenvalues can be analytically computed from the metric, damping, curvature, second covariant derivative of the potential energy, and first covariant derivative of the vector potential, a result which applies to both classical and relativistic systems. Incorporating nonlinear inequality constraints is also discussed. All derivations are tensor-based, and the computed eigenvalues themselves are coordinate-invariant, i.e., the contraction rates are independent of the chosen coordinate system. 

Simple examples including a gravity pendulum, gradient descent  with non-convex cost, Schuler dynamics, and a two-link manipulator, illustrate that the computation of the decomposed convergence rates is straightforward. The role of inequality constraints is illustrated for a controller confined to an operational envelope. 
\end{abstract}

\section{Introduction} 

Contraction analysis~\cite{lohmiller1998contraction, Lohmiller2000MechanicalContraction, lohmiller2005contraction} establishes exponential incremental convergence of a nonlinear system, and has become a standard  for solving problems in nonlinear control and estimation, and through virtual systems in synchronization~\cite{PartialContraction}. It provides a necessary and sufficient stability condition under a Riemannian metric, with closely related works including, e.g., \cite{lewis1949metric, krasovskii1963stability, hartman1961stability, boyd1985fading, DesoerHaneda1972, Pavlov2006Uniform,sontag, FB-CTDS}. In this paper, we propose a systematic analytic method to compute a suitable metric, yielding in the process analytic coordinate-invariant rates of exponential convergence.

We consider smooth real contravariant dynamics
\begin{equation}
    \dot{x}^j = {f}^j({x^k},t), \ \ \ t \ge 0 \label{eq:f}
\end{equation}
with $N$-dimensional position $x^k$ and $N$-dimensional contravariant function $f^j$. Its variational dynamics is
\begin{equation}
\delta \dot{x}^j = \frac{\partial{f^j}} {\partial x^k} \ \delta x^k \label{eq:deltadotx}
\end{equation}
We use for convenience standard tensor notation and Einstein's sum convention, which implies a summation over each index pair. Upper indices correspond to contravariant coordinates and lower indices correspond to covariant coordinates. The formulation of contraction analysis in terms of key tensor concepts, such as the covariant derivative, is summarized in Section \ref{contraction}. 

Contraction analysis provides stability conditions using a suitable metric. For a general nonlinear system, finding such a metric $M_{jk}(x^l, t)$ is not systematic, beyond finding a numerical solution to the linear matrix inequality~\cite{lohmiller1998contraction} defining the metric. This paper resolves this issue by rewriting (\ref{eq:f}) in covariant form as a natural gradient,
\begin{equation}
    M_{jk} \ \dot{x}^k \ =  \frac{\partial H}{\partial x^j} \label{eq:Mdotx}
\end{equation}
where $M_{jk}(x^h)$ and $H(x^j, t)$ may be complex, and $M_{jk}$ is symmetric and invertible. This decomposition is possible under mild integrability conditions of Frobenius \cite{bosch1987note, Frobenius1910, Lovelock} with respect to the complex linear p.d.e. $ \ \frac{\partial (M_{jk}  f^k)}{\partial x^h} = \frac{\partial (M_{hk}  f^k)}{\partial x^j} \ $ in $M_{jk} \ $. Hence the search of a contraction metric reduces to writing the system dynamics in complex natural gradient form (\ref{eq:Mdotx}).
For intuition, consider an arbitrary, real, smooth dynamics,
\begin{equation}
\dot{x}^j = f^j(x^k) \nonumber
\end{equation}
In general, this dynamics cannot be rewritten as a gradient or a natural gradient (\ref{eq:Mdotx}) using a symmetric positive definite metric $M_{jk}\ $  as in \cite{amari1998natural}, since this would imply that the cost function $H(x^k)$ always decreases and thus is a Lyapunov-like function~\cite{slotine1991applied, BeyondConvexity, Barta2022GradientSystem,EffectiveLearning2025},
\begin{equation}
\frac{d}{dt}\ H \ = \  \frac{\partial H}{\partial x^j}  M^{jk} \frac{\partial H}{\partial x^k} \ \le \ 0 \nonumber
\end{equation}
However if we allow the symmetric invertible metric to be complex and not necessarily positive definite, this can be done for any dynamics, and it also extends to non-autonomous systems. 

Section \ref{convpar} then shows that the eigenvectors of the second covariant derivative of $H(x^l, t)$ lead to an exact decomposition of the variational dynamics (\ref{eq:deltadotx}). We use the local quadratic geodesic coordinates, first introduced by C. F. Gauss~\cite{Gauss1827}, to remove a possibly coupling in form of the Christoffel term  \cite{Christoffel1869}. The computed contraction rates are the exact decoupled exponential convergence rates of (\ref{eq:deltadotx}), without any conservatism. 

Furthermore, both classical and relativistic physics \cite{Einstein:1905, einstein1915field, Hamilton} are based on a Hamiltonian function $H$, which we will show can written in the form (\ref{eq:Mdotx}). This fundamental principle of physics derives from the introduction of a differential length in a metric \cite{Einstein:1905, einstein1915field, Hamilton, Lovelock}, which implicitly defines the kinetic energy of the system. So far, little has been known about this fundamental Hamiltonian differential length being stable, indifferent or unstable over time. 
Section \ref{Ham} shows that differential lengths of general Hamiltonian dynamics correspond to exact complex analytic exponential functions, whose eigenvalues can be analytically computed from the metric, damping, curvature, second covariant derivative of the potential energy and first covariant derivative of the vector potential. This result applies to both classical and relativistic systems.


Finally, Section \ref{constrained} extends the gradient dynamics (\ref{eq:Mdotx}) to the case when positions $x^k \in  \mathbb{R}^N$ are confined to a subset $\mathbb{G}^n \subset \mathbb{R}^N$. Dirac constraint forces at the border $\partial \mathbb{G}^n$ of the set $\mathbb{G}^n$ ensure that $\mathbb{G}^n$ is invariant \cite{lohmiller2025computingquantumwavesexactly}.

 

Concluding remarks are offered in Section \ref{summary}.

\section{Contraction analysis and covariant derivatives} \label{contraction}

Contraction analysis \cite{lohmiller1998contraction}, and chaos theory \cite{lorenz1963, abarbanel1996} make extensive use of  variational displacements (\ref{eq:deltadotx}), which are differential displacements at fixed time borrowed from mathematical physics and optimization theory \cite{brysonho1975}. Whereas chaos theory and linear theory 
compute {\it numerically} the transition matrix and hence the time-averaged convergence rates in form of the Lyapunov exponents, contraction analysis provides explicit {\it analytical bounds} on the instantaneous convergence or contraction rate. 

To reach a necessary and sufficient stability condition a coordinate transformations of the form 
$$
\delta z^l = \Theta^l_k (x^j, t) \delta x^k
$$ 
has to be performed on the variational displacement similarly to a modal decomposition of LTI systems \cite{lohmiller1998contraction}. This local coordinate transformation is much more general than a coordinate change, since an explicit $z^l(x^j, t)$ does not need to exist. It generalizes the variational dynamics (\ref{eq:deltadotx}) to 
\begin{equation}
    \frac{d}{dt} \delta z^l \ = \ F_{m}^l  \delta z^m\ \ \ \ \ {\rm where} \ \ F_m^l \ = \ \left(\dot{\Theta}^l_k + \Theta^l_j \frac{\partial f^j} {\partial x^k} \right)({\Theta}_k^m)^{-1} \label{eq:Fml}
\end{equation}
The matrix $F_m^l$ is the generalized (or Riemannian) Jacobian. The corresponding differential squared length is given by
\begin{equation}
\delta z^{l} \delta z^l  = \delta x^{j} M_{jk} \delta x^k \label{eq:length}
\end{equation}
In a general curvilinear coordinate system $x^k$ the symmetric metric tensor $M_{jk}(x^l) = \Theta^l_k \Theta^l_j$ is needed to weigh the differential coordinates in the computation of the distance between two neighboring points \cite{Lovelock}. In classical physics the symmetric metric tensor $M^{jk}$ is always real and positive definite, so that any length fulfills $\delta x^{j} M_{jk} \delta x^k \ge 0$. In special or general relativity \cite{Einstein:1905} this is not the case anymore, e.g., for a symmetric but not definite Minkowski metric~\cite{Minkowski1908}. The differential length $\delta x^{j} M_{jk} \delta x^k$ can then be positive or negative with respect to the light cone \cite{Einstein:1905}. The definition above encompasses both classical and relativistic physics. It also extends to quantum physics using the conversion of action to waves in \cite{lohmiller2025computingquantumwavesexactly}.

The change of the curvilinear coordinate system $\Theta^l_j(x^h(t))$ along the velocity $\dot{x}^i$ can be computed from the metric with
\begin{equation}
    \dot{\Theta}^l_j \ = \ {\Theta}^l_k \ {\bf \gamma}_{hj}^k \ \dot{x}^h \ \ \ \ \ \ \ \ \ \ M_{jk}(x_o^l(0)) = \Theta^l_k  \Theta^l_j(x_o^l(0)) 
\end{equation}
and the Christoffel term $\gamma_{ij}^k(x^l)$ from Definition \ref{def:CoverDer} below, where $\frac{d}{dt}(\Theta^l_k \Theta^l_j) = \dot{M}_{jk}$ implies that $M_{jk} = \Theta^l_k  \Theta^l_j$ holds $\forall t \ge 0$. Summarizing the above leads to the following definition~\cite{RicciCurbastroLeviCivita1901, Lovelock}: 
\begin{definition}
\itshape
The covariant derivative of a contravariant vector $f^j(x^k,t)$ and a covariant vector $f_j = M_{jh}f^h$ are defined as
\begin{eqnarray}
   (\Theta^l_j)^{-1} \ F_m^l \ {\Theta}_k^m \ = \  f^j_{| k} &=& \frac{\partial f^j}{\partial x^k} + \gamma_{hk}^j f^h \nonumber \\
  \ M_{jh} \ f^h_{| k} \ = \    f_{j | k} &=& \frac{\partial f_j}{\partial x^k} - \gamma_{jk}^h f_h   \nonumber 
\end{eqnarray}
with the generalized Jacobian (\ref{eq:Fml}) and the Christoffel term \cite{Christoffel1869} 
\begin{equation}
\gamma_{jk}^m M_{lm} = \frac{1}{2} \left(\frac{\partial M_{jl}}{\partial x^k} + \frac{\partial M_{kl}}{\partial x^j} - \frac{\partial M_{jk}}{\partial x^l} \right) \label{eq:Christoffel}
\end{equation}
for a complex, symmetric, and invertible metric $M^{jk}(x^l) = M_{jk}^{-1}(x^l)$. 
\label{def:CoverDer}
\end{definition}
The index $_{|k}$ defines that this term is a covariant derivative with respect to the coordinate $x^k$. The covariant derivative and the generalized Jacobian (\ref{eq:Fml}) are fully equivalent. However the usage of a covariant derivative is mathematically more convenient since it can be expressed by a metric tensor only. Also recall the following linear transformation properties of tensors \cite{Christoffel1869, RicciCurbastroLeviCivita1901, Lovelock}.
\begin{lemma}
\itshape 
The covariant derivatives in Definition \ref{def:CoverDer} are tensors, since they transform linearly for any coordinate transformation $x^j(\bar{x}^l)$ as
\begin{eqnarray}
   \bar{M}_{lm} &=& \frac{\partial x^j} {\partial \bar{x}^l} \ M_{jk} \ \frac{\partial {x}^k}{\partial \bar{x}^m} \nonumber  \\ 
\bar{f}_{l | m}  &=& \frac{\partial x^j} {\partial \bar{x}^l} \ f_{j | k} \ \frac{\partial {x}^k}{\partial \bar{x}^m} \nonumber  \\  
\bar{f}^l_{| m}  &=& \frac{\partial \bar{x}^l} {\partial x^j} \ f^j_{| k} \ \frac{\partial {x}^k}{\partial \bar{x}^m} \nonumber  
\end{eqnarray}
By contrast, the Christoffel term (\ref{eq:Christoffel}) is not a tensor and transforms as \cite{Lovelock}
\begin{equation}  
     \frac{\partial {x}^n}{\partial \bar{x}^h} \bar{\gamma}_{lm}^h = \gamma_{jk}^n \frac{\partial {x}^j}{\partial \bar{x}^l} \frac{\partial {x}^k}{\partial \bar{x}^m} + \frac{\partial^2 {x}^n}{\partial \bar{x}^l \partial \bar{x}^m} \label{eq:ChristoffelTransformation}
\end{equation}
or in the standard form of contraction analysis \cite{lohmiller1998contraction}
\begin{equation}  
     \bar{\gamma}_{lm}^h {\bf T}_j^l \ \dot{x}^j \ = \ ({\bf T}_n^h)^{-1}\ \gamma_{jk}^n \ \dot{x}^j \ {\bf T}_m^k \ + \ ({\bf T}_n^h)^{-1}\ \dot{\bf T}_l^n  \label{eq:ChristoffelTransformationTime}
\end{equation}
 \label{lem:tensor}
\end{lemma}
 Note that (\ref{eq:ChristoffelTransformation}) only holds for a time-independent coordinate transformation $\bar{x}^l(x^j)$. A time-dependence can be introduced by relativistically augmenting the position $x^j \rightarrow (x^j, t)$ and the velocity $\dot{x}^j \rightarrow (\dot{x}^j, t)$ with time \cite{Einstein:1905}. This enables to introduce a pure time-varying coordine transformation $x^j - x^j_p(t)= {\bf T}_l^j(x_p^k(t), t) \bar{x}^l$ in the neighborhood of a path $x^j_p(t)$. Multiplying the Christoffel term  (\ref{eq:ChristoffelTransformation}) with ${\bf T}_j^l \dot{x}^j$  leads to (\ref{eq:ChristoffelTransformationTime})
. 


\section{Natural gradient dynamics decomposition} \label{convpar}


This section first performs an eigenvector decomposition of the covariant natural gradient dynamics (\ref{eq:Mdotx}) before the remaining Christoffel (\ref{eq:Christoffel}) coupling term is removed in local geodesic coordinate. 

Let us first use Definition~\ref{def:CoverDer} to rewrite the variational dynamics (\ref{eq:deltadotx}) of (\ref{eq:Mdotx}) as
\begin{eqnarray}
f^j_{| k}  &=& M^{hj} H_{| h k} \nonumber \\
M_{jk} (\delta \dot{x}^k + \gamma_{hi}^k \dot{x}^h \delta x^i ) &=& H_{| jk} \delta {x}^k \label{eq:covariantdelta}
\end{eqnarray}
We now recall the general definition of the characteristic equation \cite{cauchy1840memoire, Cayley1858, Frobenius1878, Lovelock}:
\begin{lemma} 
The coordinate-invariant eigenvalue matrix $\Lambda_{lm}(x^n)  = diag(\lambda_1, ..., \lambda_N)$ of a complex symmetric tensor $H_{| j k}$ with respect to a complex symmetric metric $M_{jk}$ is defined by
\begin{equation} 
{\rm det}\ (\lambda_l M_{jk}  - H_{| j k}) \ = \ 0 \nonumber
\end{equation} 
The corresponding complex orthonormal eigenvector matrix ${\bf T}_l^k(x^n)  = (T_1, ..., T_N)$ is defined with the Kronecker delta $\delta_{lm}$ \cite{kronecker1868} by
\begin{eqnarray}
    {\bf T}_m^j M_{jk} {\bf T}_l^k &=& \delta_{lm} \nonumber \\
     H_{| j k} {\bf T}_m^k &=&   M_{jk} {\bf T}_l^k \Lambda_{lm} \nonumber 
\end{eqnarray}
\label{lem:rootmetric}
\end{lemma}
In $ \ \delta \bar{x}^m = ({\bf T}_m^k)^{-1} \ \delta {x}^k \ $ coordinates, we can rewrite (\ref{eq:covariantdelta}) in the differential neighborhood of the path $x^j(t)$ as
\begin{eqnarray}
\delta_{lm} (\delta \dot{\bar{x}}^m +  \bar{\gamma}_{hi}^m \ {\bf T}_j^h \ \dot{x}^j \  \delta \bar{x}^i ) &=& \Lambda_{lm} \delta {\bar{x}}^m \label{eq:covariantdeltabarx}
\end{eqnarray}
using $ \bar{\gamma}_{hi}^m \ {\bf T}_j^h \dot{x}^j$ from Lemma \ref{lem:tensor}. This dynamics is fully decoupled up to the Christoffel term.

Now the Christoffel term of (\ref{eq:covariantdeltabarx}) can be  removed in local {\it quadratic} geodesic coordinates:
\begin{lemma}
\itshape 
In the neighborhood $\delta {x}^j = {x}^j - {x}_p^j(t)$ of the path $x_p^j(t)$, the local quadratic geodesic coordinates $\delta z^m$ are defined as
\begin{equation}
    \delta \bar{x}^l = \delta z^l - {\gamma}_{jk}^l({x}^j_p) \ \delta z^j \  \delta z^k  \nonumber
\end{equation}
with ${\gamma}_{jk}^l$ from Definition \ref{def:CoverDer}. In the local geodesic coordinates $\delta z^j$, we have in the neighborhood of the path ${x}_p^j(t)$
\begin{eqnarray}
    \bar{\gamma}_{lm}^n &=& 0 \nonumber \\
    \bar{M}_{lm}  &=& M_{lm} \nonumber 
\end{eqnarray}
\label{lem:geodesic}
\end{lemma}
\noindent {\bf Proof} \ \ C. F. Gauss ~\cite{Gauss1827, Lovelock} used Lemma \ref{lem:tensor} to perform the coordinate transformation of the Christoffel term. The metric does not change since the linear part of the geodesic coordinates is unchanged.
$ \hfill \square$

Finally Lemma \ref{lem:geodesic} allows to fully decouple (\ref{eq:covariantdeltabarx}) 
\begin{equation}
\delta_{lm} \ \frac{d}{dt}  \ \delta z^l  \ = \  \Lambda_{lm}\ \delta z^m \nonumber
\end{equation}
As a result, the decoupled exponential convergence rates are the eigenvalues $\lambda_l$ in Lemma \ref{lem:rootmetric}. This is not surprising, since in contrast to the Christoffel term, the $\lambda_l$'s are actually coordinate invariant tensors \cite{Lovelock}. Note that this decoupling, including the removal of the Christoffel term, results from the local geodesic coordinates $\delta z^m$ of Lemma \ref{lem:geodesic} being quadratic. Using linear local coordinate transformations as in the original contraction analysis \cite{lohmiller1998contraction, FB-CTDS} does not allow to remove the Christoffel term in the general case. 

\begin{figure}
    \centering
    \includegraphics[scale=0.1]{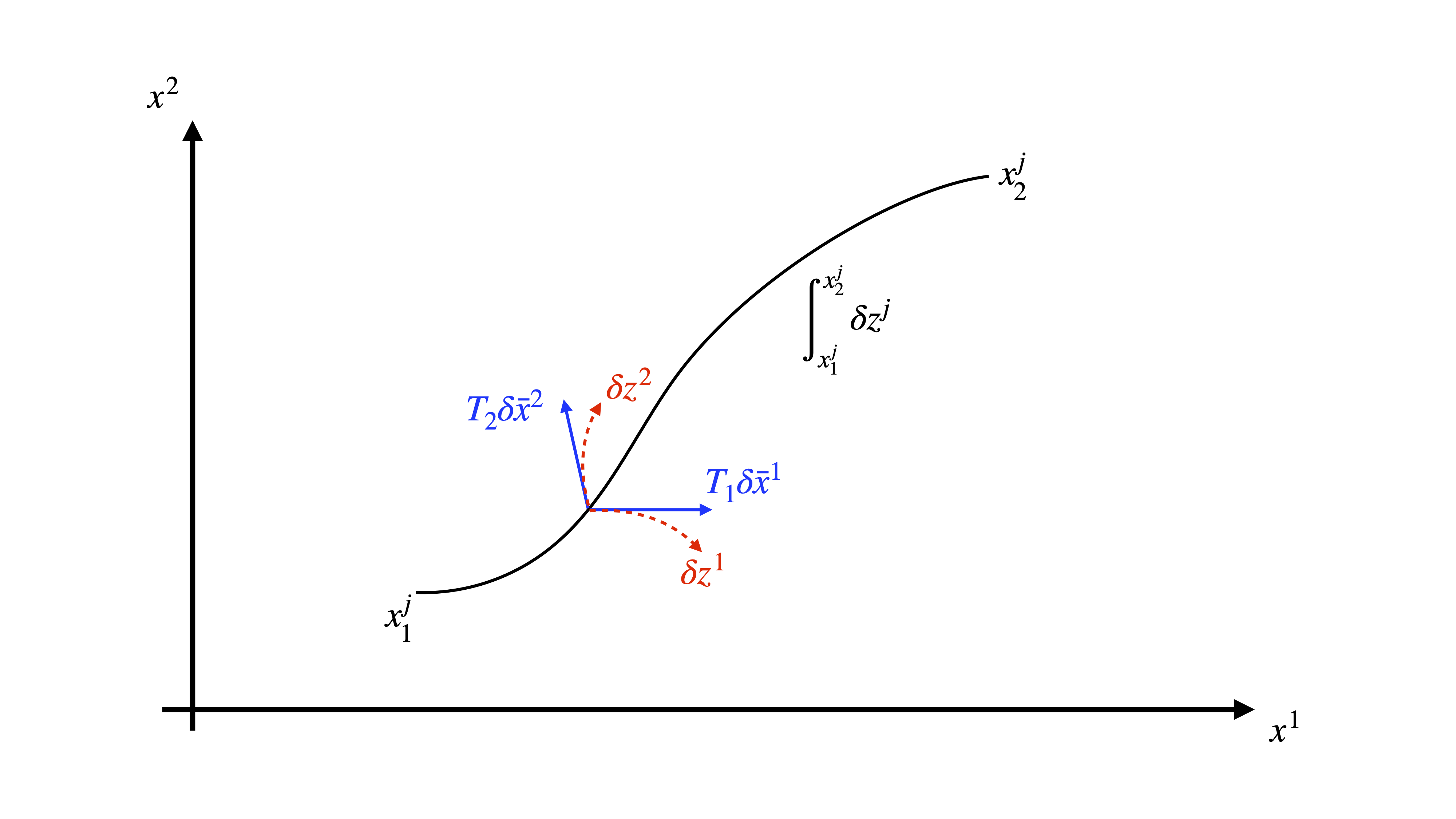}
    \caption{Eigencoordinates $\delta \bar{x}^l$, geodesic coordinates $\delta z^j$, and exponential solution between two paths $x_1^j(t), x_2^j(t)$}
    \label{fig:dxdbar}
\end{figure}

Let us summarize the above.
\begin{theorem}
The real contravariant system dynamics
\begin{equation}
   \dot{x}^j = {f}^j({x^k},t), \ \ \ x^j \in \mathbb{R}^N  \nonumber
\end{equation}
can be rewritten as covariant natural gradient dynamics
\begin{equation}
   M_{jk} \ \dot{x}^k \ = \ \frac{\partial H}{\partial x^j} \label{eq:gradientth}
\end{equation}
where the complex potential $H({x^j},t)$ and complex, symmetric and invertible metric $M_{jk}({x^n})$ are given by the decomposition
\begin{equation}
H_{| j k}(x^j)  = M_{lh} f^h_{| k} \label{eq:Fjk}
\end{equation}
The exact and decoupled contraction rates $\Lambda_{lm}(x^j)  = diag(\lambda_1, ..., \lambda_N)$ and eigenvectors ${\bf T}_l^j(x^j) = (T_1, ..., T_N)$ are given by the generalized eigenvalue equation
\begin{equation}
H_{| j k} {\bf T}_m^k =   M_{jk} {\bf T}_l^k \Lambda_{lm} \label{eq:eigenvaluesth}
\end{equation}
which imply the local geodesic coordinates $\delta z^j$ 
\begin{equation}
    \delta \bar{x}^l = ({\bf T}_l^k)^{-1}  \delta {x}^k  = \delta z^l - \bar{\gamma}_{jk}^l({x}^j) \ \delta z^j \  \delta z^k  \label{eq:localgeodesiccord}
\end{equation}
where $ \bar{\gamma}_{jk}^l({x}^j)$ from Lemma \ref{lem:tensor} is defined in $\delta \bar{x}^l$ coordinates along the path ${x}^j$.

The path integral between two arbitrary trajectories $x_1^j$ and $x_2^j$ with initial positions $x_{1o}^j$ and $x_{2o}^j$ has the exact exponential solution
\begin{equation}
    \int_{x_{1}^j}^{x_2^j} \delta z^j  = \int_{x_{1o}^j}^{x_{2o}^j}  e^{\int_o^t \lambda_j dt} \delta z^j_o \label{eq:relvector}
\end{equation}
\label{th:GenCont}
\end{theorem}
Figure \ref{fig:dxdbar} illustrates for $2$ dimensions how $(x^1, x^2)$ is decomposed in the linear eigendirections $\delta {x}^k = T_1 \delta \bar{x}^1, T_2 \delta \bar{x}^2 $ of (\ref{eq:eigenvaluesth}) and in the quadratic coordinates $\delta z^1, \delta z^1$ of (\ref{eq:localgeodesiccord}). It also illustrates the exact exponential solution (\ref{eq:relvector}) between two arbitrary paths $x_1(t)$ and $x_2(t)$.

 In the original contraction analysis \cite{lohmiller1998contraction, krasovskii1963stability, hartman1961stability, FB-CTDS} or Lyapunov approaches \cite{Lyapunov1992, slotine1991applied}, stability calculations are in general less straightforward since a metric or Lyapunov function has to be found first. The following examples illustrate how Theorem \ref{th:GenCont} provides a simple stability calculation. 
\Example{}{Consider a pendulum in a gravity field,
\begin{equation}
\left( \begin{array}{cc} 0 & i  \\ i & 0 \end{array} \right) \left( \begin{array}{c} \dot{v} \\ \dot{x} \end{array} \right)  =    
\left( \begin{array}{c} \frac{\partial H}{\partial v} \\ \frac{\partial H}{\partial x}  \end{array} \right) 
\nonumber
\end{equation}
in covariant form (\ref{eq:gradientth}) with potential $H = v^2 / 2 - \cos(x)$. 

The eigenvalue equation (\ref{eq:eigenvaluesth}) of Theorem \ref{th:GenCont} 
\begin{equation}
 \left( \left( \begin{array}{cc} 1 & 0  \\ 0 &  \cos(x) \end{array} \right) - \lambda_l   \left( \begin{array}{cc} 0 & i  \\ i & 0 \end{array} \right) \right)  T_l  = 0 \nonumber
\end{equation}
implies the decoupled contraction rates $\lambda_l = \pm \sqrt{-\cos{(x)}}$ and the eigendirections $T_l = (\pm \sqrt{-\cos{(x)}}, 1)^T$. Thus the lower equilibrium point $x=0$ is indifferent and the upper equilibrium point $x=\pi$ is unstable.}{sin} 

\Example{}{Consider the gradient descent dynamics
\begin{equation}
M_{jk} \ \dot{x}^j  = \frac{\partial H}{\partial x^k} 
\nonumber
\end{equation}
in covariant form (\ref{eq:gradientth}) with cost function $H(x^k,t)$ and metric $M_{jk}(x^l)$.

The eigenvalue equation (\ref{eq:eigenvaluesth}) of Theorem \ref{th:GenCont} 
\begin{equation}
\left( H_{| j k} + \lambda_l  M_{jk} \right) T_l   = 0 \nonumber
\end{equation}
implies the decoupled contraction rates $\lambda_l$. 

Note that the original contraction analysis \cite{lohmiller1998contraction} has to conservatively take $\max(\lambda_l)$ as contraction rate, since the variational dynamics is not decoupled in \cite{lohmiller1998contraction}. According to (\ref{eq:relvector}) in Theorem \ref{th:GenCont}, neighboring trajectories converge exactly with $e^{\int_o^t \lambda_l dt}$.}{gradient} 

Note that all lemmas in \cite{lohmiller1998contraction} on adaptive control, observer design, or combination principles can be applied to Theorem \ref{th:GenCont} as well.

\section{Exact decomposition of Hamiltonian dynamics} \label{Ham}

This section applies the covariant natural gradient (\ref{eq:gradientth}) of Theorem \ref{th:GenCont} to classical or relativistic Hamiltonian dynamics \cite{Hamilton, Einstein:1905, einstein1915field, Goldstein} 
\begin{eqnarray}
  \left( \begin{array}{cc} 0 & i \delta_j^k  \\ i \delta_k^j  &  D_{jk} \end{array} \right)
\left( \begin{array}{c} \dot{p}_j \\ \dot{q}^j \end{array} \right)  &=&    \left( \begin{array}{c} \frac{\partial H}{\partial p_k}  \\ \frac{\partial H}{\partial q^k} \end{array} \right), \ \ \ t \ge 0 \label{eq:Ham} 
\end{eqnarray}
with Hamiltonian 
\begin{equation}
    H = \frac{1}{2} (p_j - Q \ A_j) M^{jk}(q^l) (p_k - Q \ A_k) + V(q^l,t) \label{eq:Hamiltonian}
\end{equation}
with classical or relativistic {\it real} invertible inertia tensor $M^{jk}(q^l)$ \cite{Minkowski1908, einstein1915field, Lovelock}, {\it real} potential energy $V(q^l,t)$, {\it imaginary} vector potential $A_k(q^l, t)$, {\it real} constant charge $Q$, $N$-dimensional {\it real} position $q^k$ and {\it imaginary} momentum $p_k$, Kronecker delta $\delta_k^j = M_{kh} M^{jh}$ \cite{kronecker1868}, and constant {\it real} damping matrix $D_{jk}(q^l,t)$ with $D_{jk|l}=0$. Note that while the asssociated Lagrangian second-order dynamics in $q^k$ is fully real, for convenience the momentum $p_k$ and vector potential $A_k$ are defined above as imaginary to later allow for a complex eigenvector computation. The dynamics (\ref{eq:Ham}) can be rewritten as
\begin{eqnarray}
 \left( \begin{array}{cc} 0 & i \delta_j^k  \\ i \delta_k^j  &  D_{jk} \end{array} \right)
\left( \begin{array}{c} \dot{p}_j - \gamma_{jh}^l  \dot{q}^h p_l \\ \dot{q}^j \end{array} \right)  &=&    \left( \begin{array}{c} M^{jk} (p_k- Q A_k)  \\ \frac{\partial V}{\partial q_k} - {p}_j Q A_{|k}^j  \end{array} \right) \label{eq:Ham2}
\end{eqnarray}
Taking the variation of (\ref{eq:Ham2}) yields \cite{lohmiller2012exact} 
\begin{equation}
\left( \begin{array}{cc} 0 & i \delta_j^k  \\ i \delta_k^j  &  D_{jk} \end{array} \right)
    \  \left( \begin{array}{c} \frac{d}{dt} \ D p_j- \gamma_{jh}^l \dot{q}^h D p_l \\ \delta \dot{q}^j + \gamma_{hk}^j \dot{q}^h \delta q^k \end{array} \right)  =   \left( \begin{array}{cc} M^{jk} & - Q A_{|k}^j \\ -Q A_{|j}^k & V_{| j k} - R_{j \ hk}^{\ l} p_l \dot{q}^h \end{array} \right)  \left( \begin{array}{c} D p_k \\ \delta q^k  \end{array} \right) \label{eq:VarHam}
\end{equation}
with $D {p}_j = p_{j|k} \delta{q}^k = \delta {p}_j - \gamma_{jk}^l p_l \delta{q}^k$ and where we used \cite{lohmiller2012exact}
\begin{eqnarray}
     \frac{d}{dt} \ D p_j- \gamma_{jh}^m \dot{q}^h D p_m &=& \frac{d}{dt} (\delta {p}_j - \gamma_{jk}^l p_l \delta{q}^k) - \gamma_{jh}^m \dot{q}^h (\delta {p}_m - \gamma_{km}^l p_l  \delta{q}^k)  \nonumber \\
    &=&  - \left( i V_{| j k} + R_{jlkh} \dot{q}^l \dot{q}^h \right) \delta q^k - Q A_{|j}^k \delta p_k \nonumber
\end{eqnarray}
with the  Riemann curvature tensor~\cite{Riemann1867, Lovelock}
\begin{equation}
    R_{jlkh} = \left( \frac{\partial \gamma_{jk}^m}{\partial q^h} - \frac{\partial \gamma_{jh}^m}{\partial q^k} + 
\gamma_{jk}^m \gamma_{mh}^m  -  \gamma_{jh}^m  \gamma_{km}^m \right) H_{ml}  \label{eq:curvature}
\end{equation} 
Since (\ref{eq:VarHam}) has the same form as (\ref{eq:covariantdelta}), Theorem \ref{th:GenCont} yields the following result.
\begin{theorem}
Consider the Hamiltonian dynamics in covariant natural gradient form
\begin{eqnarray}
  \left( \begin{array}{cc} 0 & i \delta_j^k  \\ i \delta_k^j  &  D_{jk} \end{array} \right)
\left( \begin{array}{c} \dot{p}_j \\ \dot{q}^j \end{array} \right)  &=&    \left( \begin{array}{c} \frac{\partial H}{\partial p_k}  \\ \frac{\partial H}{\partial q^k} \end{array} \right) , \ \ \ t \ge 0 \label{eq:Hamth} 
\end{eqnarray}
with Hamiltonian (\ref{eq:Hamiltonian}), Kronecker delta $\delta_k^j = M_{kh} M^{jh}$ \cite{kronecker1868}, constant {\it real} damping matrix $D_{jk}(q^l,t)$ with $D_{jk|l}=0$ and augmented state vector ${x}^k = ( p_k, q^k)$.

The exact and decoupled contraction rates $\Lambda_{lm}(x^j)  = diag(\lambda_1, ..., \lambda_N)$ and eigenvectors ${\bf T}_l^j(x^j) = (T_1, ..., T_N)$ of Theorem \ref{th:GenCont}  are given by the generalized eigenvalue equation
\begin{equation}
   \left( \begin{array}{cc} M^{jk} & - Q A_{|k}^j \\ -Q A_{|j}^k & V_{| j k} - R_{jlkh} \dot{q}^l \dot{q}^h  \end{array} \right) {\bf T}_m^j  =    \left( \begin{array}{cc} 0 & i \delta_j^k  \\ i \delta_k^j  & D_{jk} \end{array} \right)  {\bf T}_l^j \Lambda_{lm} \label{eq:eigenHamth}
\end{equation}
\label{th:HamCont}
\end{theorem}
Theorem \ref{th:HamCont} generalizes the combined Lyapunov and contraction analysis of Hamiltonian dynamics in \cite{lohmiller2012exact} to a pure and fully decoupled contraction analysis. In contrast to the energy based Lyapunov proofs \cite{Lyapunov1992, slotine1991applied} Theorem \ref{th:HamCont} exactly decomposes the dynamics in its eigendynamics which allows to compute an exact exponential solution, even in the nonlinear and non-autonomous case. 

\begin{figure}
    \centering
    \includegraphics[scale=0.5]{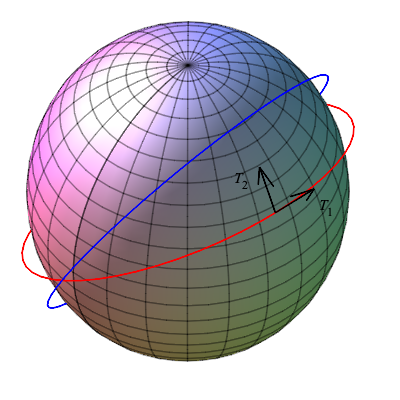}
    \caption{Two orbiting satellites with eigendirections $T_1, T_2$ }
    \label{fig:satellte}
\end{figure}

\Example{}{Consider the spherical motion of a satellite
\begin{equation}
 \left( \begin{array}{cc} 0 &  i \delta_j^k  \\  i \delta_j^k  & 0 \end{array} \right) \left( \begin{array}{c} \dot{p}_j \\ \dot{q}^j \end{array} \right)  =  \left( \begin{array}{c} \frac{\partial H}{\partial p_k}  \\ \frac{\partial H}{\partial q^k} \end{array} \right) \label{eq:satmotion}
\end{equation}
with Hamiltonian $H = -\frac{1}{2 M} p_j M^{jk}(q^l) p_k$, position $q^j = (\phi, \psi)$, latitude $\phi$, longitude $\psi$, imaginary momentum $p_l = i M^{lh} \dot{q}^h$, satellite mass $M$ and metric tensor with Christoffel term
$$
M_{jk} = R^2 \left(
\begin {array}{cc} 
1 & 0 \\
0 & \cos{\phi}^2
\end {array} 
\right), \ \ \ 
\gamma_{jk}^1  =
\left(
\begin{array}{cc}
0 & 0\\
0 & \sin \phi \cos \phi
\end{array}
\right), 
\gamma_{jk}^2  =
\left(
\begin{array}{cc}
0 & -\tan \phi \\
-\tan \phi & 0
\end{array}
\right)
$$
We can compute e.g. with \cite{christoffel_symbols_calculator} the curvature tensor (\ref{eq:curvature})
$$
R_{jlkh} \dot{q}^l \dot{q}^h = R^2 \cos^2 \phi \left(
\begin {array}{cc} 
\dot{\psi}^2 & -\dot{\psi} \dot{\phi} \\
-\dot{\psi} \dot{\phi} & \dot{\phi}^2
\end {array} 
\right) 
$$ 
The contraction rates (\ref{eq:eigenHamth}) of Theorem \ref{th:HamCont} 
\begin{equation}
  \left( \left( \begin{array}{cc} M^{jk} & 0 \\ 0 & - R_{jlhk} \dot{q}^l \dot{q}^h \end{array} \right) - \lambda_l  \left( \begin{array}{cc} 0 & i \delta_j^k  \\ i \delta_k^j  & 0 \end{array} \right) \right) T_l  = 0 \nonumber
\end{equation}
are $\lambda_1 = 0$ in the flight direction $T_1 = (\dot{\phi}, \dot{\psi}\cos \phi)^T$ and $\lambda_2 = \pm i \sqrt{ (\dot{\psi}^2 \cos^2 \phi + \dot{\phi}^2) }$ orthogonal to the flight direction $T_2 = ( \dot{\psi}\cos \phi, -\dot{\phi})^T$ in Figure \ref{fig:satellte}. Figure \ref{fig:satellte} illustrates that the lateral osciallation comes from the Earth curvature and that distances in longitudinal direction are preserved.

Let us now compute the estimated position $q^j$ of a satellite with a measurement of the Cartesian gravity vector $g_j(t)$. The measurement implies the Hamiltonian dynamics (\ref{eq:satmotion}) with $H = -\frac{1}{2 M} p_j M^{jk}(q^l) p_k + V$,  potential energy $V = - g_j x^j$ and Cartesian position ${\bf x} = R ( \cos \psi \cos \phi, \sin \psi \cos \phi, \sin \phi)$. We can compute e.g. with \cite{christoffel_symbols_calculator}
$$
V_{| j k} = \frac{\partial^2 V}{\partial x^j \partial x^k} - \gamma_{jk}^l \frac{\partial V}{\partial x^l} = -\frac{V}{R^2} M_{jk}
$$
The contraction rates (\ref{eq:eigenHamth}) of Theorem \ref{th:HamCont} 
\begin{equation}
  \left( \left( \begin{array}{cc} M^{jk} & 0 \\ 0 &  V_{| j k} - R_{jlkh} \dot{q}^l \dot{q}^h  \end{array} \right) - \lambda_l  \left( \begin{array}{cc} 0 &  i \delta_j^k  \\ i \delta_k^j  & 0 \end{array} \right) \right) T_l = 0 \nonumber
\end{equation}
are now $\lambda_1 = \pm i \sqrt{\frac{V}{R^2} }$ in the  flight direction and $\lambda_2 = \pm i \sqrt{-\frac{V}{R^2} + (\dot{\psi}^2 \cos^2 \phi + \dot{\phi}^2) }$ orthogonal to the flight direction in Figure \ref{fig:satellte}.

Note that the metric above was introduced in \cite{Schmalz2007Analytic} to remove the virtual effect of Coriolis or transport forces in an earlier study \cite{titterton2004strapdown}.  However, in \cite{Schmalz2007Analytic} the real Earth curvature effect was still neglected. The above calculation is now exact in computing the Schuler frequency \cite{Schuler1923}, as it fully accounts for  the effect of the Earth curvature. }{SDAball}

\Example{}{Consider a two-link robot manipulator 
\begin{equation}
\left( \begin{array}{cc} 0 & i \delta_j^k  \\ i \delta_k^j  &  D_{jk} \end{array} \right)
\left( \begin{array}{c} \dot{p}_j \\ \dot{q}^j \end{array} \right)  =  \left( \begin{array}{c} \frac{\partial H}{\partial p_k}  \\ \frac{\partial H}{\partial q^k} \end{array} \right) \nonumber
\end{equation}
with Hamiltonian $H = \frac{1}{2} p_j M^{jk}(q^l) p_k + V$, potential energy $V(q^k, t)$, constant damping matrix $D_{jk}(q^l,t), D_{jk|l}=0$, the periodic link angle $-\pi \le q^1, q^2 \le \pi$ of link $1$ and $2$ and the related imaginary momentum $p_l = i M^{lh} \dot{q}^h$. The Cartesian end-position of link $1$ and $2$ is
\begin{equation}
{\bf x}^1 = l \left( \begin{array}{c} \cos q^1 \\ \sin q^1 
\end{array}
\right) \ \ \ \ {\bf x}^2 = l \left( \begin{array}{c} \cos q^1 + \cos
    q^2 \\ \sin q^1 + \sin q^2 \end{array} \right) \nonumber
\end{equation}
with link length $l$ and the effective mass $M$ implies the inertia tensor
\begin{eqnarray}
M_{jk} &=& m l^2 \left(
\begin{array}{cc} 2 & \cos(q^2 -q^1) \\ \cos(q^2 - q^1) & 1
\end{array} \right) \nonumber 
\end{eqnarray}
The robot kinematics corresponds to a closed $2$-dimensional surface
within a $4$-dimensional space. The curvature tensor (\ref{eq:curvature}) can be computed e.g. with \cite{christoffel_symbols_calculator} as
$$
R_{jlkh} \dot{q}^l \dot{q}^h  = m l^2
\frac{\cos(q^2 - q^1)} {2-\cos(q^2-q^1)}
\left( 
\begin {array}{cc}  
(\dot{q}^2)^2 & -\dot{q}^1 \dot{q}^2 \\
-\dot{q}^1 \dot{q}^2 & (\dot{q}^1)^2 
\end {array} \right)
$$
Hence the curvature under motion of the open-loop 2-link robot acts as an indifferent (unstable) spring for $\cos (q^2 -q^1) \ge 0$ ($\cos(q^2-q^1) \le 0$). The contraction rates (\ref{eq:eigenHamth}) of Theorem \ref{th:HamCont} are
\begin{equation}
   \left( \left( \begin{array}{cc} M^{jk} & 0 \\ 0 & V_{| j k} - R_{jlkh} \dot{q}^l \dot{q}^h   \end{array} \right) - \lambda_l  \left( \begin{array}{cc} 0 & i \delta_j^k  \\ i \delta_k^j  & D_{jk} \end{array} \right) \right) T_l = 0 \nonumber
\end{equation} 
The second covariant derivative of the potential energy $V_{| j k}$ and the damping matrix $D_{jk}$ can be used the place the contraction rates $\lambda_l$ of the closed-loop dynamics.}{robot}

\section{Exact decomposition of constrained dynamics} \label{constrained}

We now extend the covariant dynamics (\ref{eq:gradientth}) of Theorem \ref{th:GenCont} and \ref{th:HamCont} to systems with spatial inequality constraints. 
\begin{definition} 
The constrained multiply connected manifold $ \ \mathbb{G}^N \subseteq \mathbb{R}^N$ is defined by the $g=1, ..., G$ inequality constraints $ \ f_g({\bf x})  \le  0 \ $.
The set of active constraints $\mathbb{G} \subseteq  \{1, ..., G \}$  is the set of indices $g$  on the boundary $\partial \mathbb{G}^N$ of $\mathbb{G}^N$ such that $ \ f_g({\bf x}) = 0 \ $.
\label{def:constraint}
\end{definition}
The constrained dynamic (\ref{eq:gradientth}) is then of the form \cite{bryson1975applied}
\begin{eqnarray}
  {M}_{jk} \dot{x}^k &=& \frac{\partial H}{\partial x^j}  + \sum_{g \in \mathbb{G}} \frac{\partial f_g}{\partial x^j}  \mu_g = \frac{\partial H}{\partial x^j} + \sum_{g=1}^G  \frac{\partial f_g}{\partial x^j}  h(f_g) \mu_g  \label{eq:Mxdotg}
\end{eqnarray}
with the Lagrange multipliers $\mu_g \in \mathbb{R}$ \cite{Lagrange}, Heaviside step function $h$, and Dirac impulse $\delta (f_g) = \frac{\partial h(f_g)}{\partial f_g}$. The constraints $g \in \mathbb{G}$ are respected at $t + dt$ if \cite{bryson1975applied}
\begin{eqnarray}
 \dot{f}_g = \frac{\partial f_g }{\partial {x}^k} {M}^{jk}  \left( \frac{\partial H}{\partial x^j} + \sum_{g \in \mathbb{G}} \frac{\partial f_g}{\partial x^j}  \mu_g  \right)  &\le&  0 \nonumber
\end{eqnarray}
We thus have the following result.
\begin{lemma} 
The Lagrange multipliers are
\begin{eqnarray}
 \mu_g &=& - \left( \frac{\partial f_g }{\partial {x}^k} {M}^{jk} \sum_{g \in \mathbb{G}}  \frac{\partial f_g }{\partial {x}^j} \right)^{-1} (e+1)\frac{\partial f_g }{\partial {x}^k} {M}^{jk}  \frac{\partial H}{\partial x^j} \label{eq:lambda}
\end{eqnarray}
with the coefficient of restitution $e \ge 0$. We have $ \ e = 1 \ $ for a perfectly elastic collision, $ \ e = 0 \ $ for  a plastic collision, and $\ 0 < e < 1 \ $ for a partially elastic collision.
 \label{lem:activeconstraint}
\end{lemma}
Note that initial principles for introducing convex constraints in contraction analysis were illustrated in \cite{lohmiller1999nonlinear, tabareau2007geometry, nguyen2017contraction}. The total variation of the unconstrained dynamics (\ref{eq:Mxdotg}) can be written as 
\begin{equation}
{M}_{jk} \frac{d}{dt} \ \delta {x}^j =  {H}_{| j k} \ \delta {x}^j  + \sum_{g = 1}^G \left( h(f_g) \left( f_{g| j k}\lambda_k + \frac{\partial f_g}{\partial {x}^k} \frac{\partial \lambda_k }{\partial {x}^j}  \right)  + \delta(f_g) \mu_g \frac{\partial f_g}{\partial {x}^k}  \frac{\partial f_g}{\partial {x}^j} \right) \delta {x}^j \nonumber
\end{equation}
The term weighted with $h(f_g)$ can be neglected for an elastic collision since the collision time is infinitesimal small. The term weighted with $\delta(f_g)$ normalizes the variational displacement component $\frac{\partial f_g}{\partial {x}^j} \delta {x}^j$ to $0$ at a plastic collision and to $-\frac{\partial f_g}{\partial x^j} \delta {x}^j$ at a fully elastic collision (see Lemma \ref{lem:activeconstraint}). 

Summarizing the above yields a more general version of Theorem \ref{th:GenCont} and \ref{th:HamCont}.
\begin{theorem}
The real contravariant system dynamics
\begin{equation}
    \dot{x}^j = {f}^j({x^k},t), \ \ \ x^j \in \mathbb{R}^N  \nonumber
\end{equation}
can be rewritten with (partial) elastic collisions $x^j \in \mathbb{G}^n \label{eq:Mphi}$ (see Definition  \ref{def:constraint}) as constrained covariant natural gradient dynamics
\begin{equation}
    M_{jk} \ \dot{x}^k \ = \ \frac{\partial H}{\partial x^j} + \sum_{g  \in \mathbb{G}} \frac{\partial f_g}{\partial x^j}  \mu_g , \ \ \ x^j \in \mathbb{G}^N \label{eq:gradientthconst}
\end{equation}
where the complex potential $H({x^j},t)$ and complex, symmetric and invertible metric $M_{jk}({x^n})$ is given by  
(\ref{eq:Fjk}) of Theorem \ref{th:GenCont} and the set of active constraints $\mathbb{G}$ is given in Definition \ref{def:constraint}.

In a collision free path segment, the path integral between two arbitrary trajectories $x_1^k$ and $x_2^k$ exponentially converges with contraction rates (\ref{eq:eigenvaluesth}) according to Theorem \ref{th:GenCont} and \ref{th:HamCont}. At the collision time instance the variational displacement component $\frac{\partial f_g}{\partial x^j} \delta {x}^j$ is set to $0$ at a plastic collision $e=0$ and to $-\frac{\partial f_g}{\partial x^j} \delta {x}^j$ at a fully elastic collision $e=1$ in Lemma \ref{lem:activeconstraint}.
\label{th:GenConst}
\end{theorem}
Note that after a plastic collision the dimensionality $N$ of the metric $M_{j k}$ (\ref{eq:gradientth}) is reduced such that the contraction rates before and after the collision can be computed separately with two different metrics $M_{j k}$. 

\Example{}{Consider the second-order dynamics
\begin{eqnarray}
\left( \begin{array}{cc}
     0 & i \\
     i & 1
\end{array} \right) 
\left( \begin{array}{c}
     \dot{p}  \\
     \dot{q}
\end{array} \right) &=& 
\left( \begin{array}{c}
     p  \\
     q
\end{array} \right)
- \left( \begin{array}{c}
     0 \\
     u(t)
\end{array} \right) + \sum_{g \in \mathbb{G}} \frac{\partial f_g}{\partial x^j} \mu_g
\nonumber
\end{eqnarray}
limited with the active constraints $\mathbb{G} \subseteq  \{1, 2\}$ of Definition \ref{def:constraint}
\begin{equation}
    -2 \le q \le 2 \ \ \ \ \ \Leftrightarrow \ \ \ \ \ f_1 = q-2 \le 0\ , \ f_2 = -q+2 \le 0 \nonumber
\end{equation}

The contraction rate is given by the eigenvalue equation (\ref{eq:eigenvaluesth}) of Theorem \ref{th:GenCont}
$$
\left( \left( \begin{array}{cc}
     1 & 0 \\
     0 & 1
\end{array} \right) - \left( \begin{array}{cc}
     0 & i \\
     i & 1
\end{array} \right) \lambda_l \right) T_l = 0 
$$
which implies the double contraction rate $\lambda_l = - \frac{1}{2}$. According to Theorem \ref{th:GenConst}, the elastic collision is indifferent for $e=1$ at the collision time instance. Note that many practical systems have similar linear constraints, such as end positions or actuator limits.}{spring}

\section{Concluding Remarks} \label{summary}

Theorem \ref{th:GenCont} computes the exact and decomposed exponential convergence rates of general nonlinear time-varying systems based on contraction analysis \cite{lohmiller1998contraction}. A natural metric is found by writing the general system dynamics as a complex natural gradient dynamics in covariant form. Similarly, Theorem~\ref{th:HamCont} computes the analytic exponential convergence rates of the differential length of a general Hamiltonian dynamics. Both theorems correspond to an exact eigencoordinate decomposition whose variational dynamics can be solved analytically. Theorem \ref{th:GenConst} extends these results to systems with nonlinear inequality constraints, where Dirac collision forces ensure that the system remains in a bounded space. 



The computed eigendirections are tensors and the computed eigenvalues are coordinate-invariant, i.e., we get the same contraction rate in any coordinate system.  Diverse examples illustrate that the computation of the decomposed convergence rates is straightforward. Current research focuses on nonlinear optimal observer design and on the discrete-time and continuous-discrete-time cases, with applications to machine learning and closed-loop control design.  Also, the decomposition of Hamiltonian dynamics allows a decomposed action computation. These decomposed actions can then be used in quantum physics to compute the quantum wave \cite{lohmiller2025computingquantumwavesexactly}, e.g., for systems with nonlinear potentials.

\noindent {\bf Acknowledgements} \ \ This paper benefited from discussions with Airan Thallemer.

\normalem
\bibliographystyle{abbrv}
{\fontsize{10}{0} \selectfont \bibliography{References}{}}

@article{amari1998natural,
  title={Natural gradient works efficiently in learning},
  author={Amari, Shun-Ichi},
  journal={Neural Computation},
  volume={10},
  number={2},
  pages={251--276},
  year={1998},
  publisher={MIT Press}
}

@book{abarbanel1996,
  title={Analysis of Observed Chaotic Data},
  author={Abarbanel, Henry D. I.},
  isbn={0387983724},
  lccn={95050476},
  series={Institute for Nonlinear Science},
  year={1996},
  publisher={Springer-Verlag},
  address={New York, NY}
}

@article{Barta2022GradientSystem,
  title={Every ordinary differential equation with a strict Lyapunov function is a gradient system},
  author={B{\'a}rta, Tom{\'a}{\v{s}} and Chill, Ralph and Fa{\v{s}}angov{\'a}, Eva},
  journal={Journal of Differential Equations},
  year={2022}
}

@article{bosch1987note,
  author    = {A. J. Bosch},
  title     = {{Note on the factorization of a square matrix into two Hermitian or symmetric matrices}},
  journal   = {SIAM Review},
  volume    = {29},
  number    = {3},
  pages     = {463--466},
  year      = {1987}
}

@book{brysonho1975,
  author    = {Arthur E. Bryson, Jr. and Yu-Chi Ho},
  title     = {Applied Optimal Control: Optimization, Estimation, and Control},
  publisher = {Hemisphere Publishing Corporation},
  address   = {Washington, DC},
  year      = {1975},
}

@book{Pavlov2006Uniform,
  author    = {Alexey Pavlov and Nathan van de Wouw and Henk Nijmeijer},
  title     = {Uniform Output Regulation of Nonlinear Systems: A Convergent Dynamics Approach},
  publisher = {Birkh{\"a}user},
  address   = {Boston},
  year      = {2006},
  series    = {Systems \& Control: Foundations \& Applications},
  isbn      = {978-0-8176-4445-1}
}

@Book{FB-CTDS,
  author    = {F. Bullo},
  title     = {Contraction Theory for Dynamical Systems},
  year      = {2024},
  edition   = {1.2},
  publisher = {Kindle Direct Publishing},
  ISBN      = {979-8836646806},
  url       = {https://fbullo.github.io/ctds},
}

@article{boyd1985fading,
  title={Fading Memory and the Problem of Approximating Nonlinear Operators with Volterra Series},
  author={Boyd, Stephen and Chua, Leon O.},
  journal={IEEE Transactions on Circuits and Systems},
  volume={32},
  number={11},
  pages={1150--1161},
  year={1985}
}

@book{bryson1975applied,
  title = {Applied Optimal Control: Optimization, Estimation and Control},
  author = {Bryson, Arthur E. and Ho, Yu-Chi},
  year = {1975},
  publisher = {Hemisphere Publishing Corp.},
  address = {New York}
}

@misc{christoffel_symbols_calculator,
  author = {Dhananjhay Bansal},
  title = {Christoffel Symbols Calculator},
  howpublished = {\url{https://christoffel-symbols-calculator.com}},
  note = {Accessed: 2025-08-17}
}

@article{cauchy1840memoire,
  author = {Cauchy, Augustin-Louis},
  title = {Mémoire sur l'intégration des équations linéaires},
  journal = {Comptes Rendus de l'Académie des Sciences de Paris},
  year = {1840},
  volume = {11},
  pages = {1--12},
  note = {Oeuvres complètes, Série 1, Tome 8}
}

@article{Cayley1858,
  author = {Arthur Cayley},
  title = {A Memoir on the Theory of Matrices},
  journal = {Philosophical Transactions of the Royal Society of London},
  year = {1858},
  volume = {148},
  pages = {17--37},
  url = {https://en.wikipedia.org/wiki/Arthur_Cayley}
}

@article{Christoffel1869,
  author = {Christoffel, Elwin Bruno},
  title = {{Über die Transformation der homogenen Differentialausdrücke zweiten Grades}},
  journal = {{Journal für die reine und angewandte Mathematik}},
  volume = {70},
  pages = {46--70},
  year = {1869},
  url = {http://gdz.sub.uni-goettingen.de/dms/load/img/?PPN=GDZPPN002153882&IDDOC=266356},
  language = {de}
}

@article{DesoerHaneda1972,
  author = {C. A. Desoer and H. Haneda},
  title = {The measure of a matrix as a tool to analyze computer algorithms for circuit analysis},
  journal = {IEEE Transactions on Circuit Theory},
  volume = {19},
  number = {},
  pages = {480-486},
  year = {1972}
}

@article{Einstein:1905,
author = {Einstein, Albert},
title = {{Zur Elektrodynamik bewegter Körper}},
journal = {Annalen der Physik},
volume = {322},
number = {10},
pages = {891--921},
year = {1905},
doi = {10.1002/andp.19053221004},
url = {https://de.wikipedia.org/wiki/Spezielle_Relativit%C3%A4tstheorie}
}

@article{einstein1915field,
  author = {Einstein, Albert},
  title = {{Die Feldgleichungen der Gravitation}},
  journal = {Sitzungsberichte der Königlich Preußischen Akademie der Wissenschaften (Berlin)},
  year = {1915},
  pages = {844--847},
  url = {https://einsteinpapers.press.princeton.edu/vol6-doc21}
}

@article{Frobenius1910,
author = {Ferdinand Georg Frobenius},
title = {{Über die mit einer Matrix vertauschbaren Matrizen}},
journal = {Sitzungsberichte der Preussischen Akademie der Wissenschaften},
year = {1910}
}

@article{Frobenius1878,
  author    = {Ferdinand Georg Frobenius},
  title     = {{Über lineare Substitutionen und bilineare Formen}},
  journal   = {Journal für die reine und angewandte Mathematik},
  year      = {1878},
  volume    = {84},
  pages     = {1-63},
}

@article{Gauss1827,
  author    = {Carl Friedrich Gauss},
  title     = {Disquisitiones Generales Circa Superficies Curvas},
  journal   = {Commentationes Societatis Regiae Scientiarum Gottingensis Recentiores},
  volume    = {6},
  pages     = {99--146},
  year      = {1827},
  url       = {https://www.gutenberg.org/files/36856/36856-pdf.pdf}
}

@book{Goldstein,
    title = {{Classical Mechanics}},
    author = {Goldstein, H.},
    year = {1980},
    publisher = {Addison-Wesley}
}

@inproceedings{Hamilton,
  title        = {{Second essay on a general method in dynamics}},
  author       = {Hamilton, R.W.},
  year         = 1835,
  booktitle    = {Philosophical Transactions of the Royal Society}
}

@article{hartman1961stability,
  author    = {Philip Hartman},
  title     = {On stability in the large for systems of ordinary differential equations},
  journal   = {Canadian Journal of Mathematics},
  volume    = {13},
  pages     = {480--492},
  year      = {1961}
}

@article{kronecker1868,
  author = {Leopold Kronecker},
  title = {{Über bilineare Formen}},
  journal = {{Journal für die reine und angewandte Mathematik}},
  volume = {68},
  pages = {273--285},
  year = {1868},
  language = {German}
}

@article{EffectiveLearning2025,
  title={Effective Learning Rules as Natural Gradient Descent},
  author={Lucas Shoji and Kenta Suzuki and Leo Kozachkov},
  journal={Neural Computation},
  year={2025},
  pages={1-26},
  doi={10.1162/NECO.a.1474}
}

@book{krasovskii1963stability,
  title={Stability of Motion},
  author={Krasovskii, N. N.},
  year={1963},
  publisher={Stanford University Press},
  note={Original Russian edition: 1959}
}

@inproceedings{lohmiller2012exact,
  author    = {Winfried Lohmiller and Jean-Jacques Slotine},
  title     = {Exact Modal Decomposition of Non-linear Second-Order Systems},
  booktitle = {AIAA Guidance, Navigation, and Control Conference},
  year      = {2012},
  doi       = {10.2514/6.2012-4751}
}

@article{lohmiller2005contraction,
  title = {Contraction Analysis of Nonlinear Distributed Systems},
  author = {Lohmiller, Winfried and Slotine, Jean-Jacques},
  journal = {International Journal of Control},
  year = {2005},
  volume = {78},
  number = {9},
  pages = {678--688},
  doi = {10.1080/00207170500130952}
}

@article{Lohmiller2000MechanicalContraction,
  author  = {Winfried Lohmiller and Jean-Jacques Slotine},
  title   = {Control System Design for Mechanical Systems Using Contraction Theory},
  journal = {IEEE Transactions on Automatic Control},
  year    = {2000},
  volume  = {45},
  number  = {5},
  pages   = {884--889}
}

@book{Lagrange,
    title = {{Mécanique analytique}},
    author = {Lagrange, J.L.},
    year = {1788},
    publisher = {Chez la veuve Desaint à Paris}
}

@article{lewis1949metric,
  author    = {D. C. Lewis},
  title     = {Metric properties of differential equations},
  journal   = {American Journal of Mathematics},
  volume    = {71},
  pages     = {294--312},
  year      = {1949}
}

@article{lohmiller1998contraction,
  title={On contraction analysis for nonlinear systems},
  author={Lohmiller, Winfried and Slotine, Jean-Jacques},
  journal={Automatica},
  volume={34},
  number={6},
  year={1998},
  publisher={Elsevier},
  doi={10.1016/S0005-1098(98)00019-3}
}

@article{lohmiller1999nonlinear,
  title={Nonlinear process control using contraction theory},
  author={Lohmiller, Wolfgang and Slotine, Jean-Jacques},
  journal={AIChE Journal},
  volume={46},
  number={3},
  pages={588--596},
  year={1999},
  doi={10.1002/aic.690460317}
}

@misc{lohmiller2025computingquantumwavesexactly,
      title={On computing quantum waves exactly from classical and relativistic action}, 
      author={Winfried Lohmiller and Jean-Jacques Slotine},
      year={2026},
      eprint={2405.06328},
      archivePrefix={arXiv},
      primaryClass={quant-ph},
      url={https://arxiv.org/abs/2405.06328}, 
}

@book{Lovelock,
    title = {Tensors, Differential Forms, and Variational Principles},
    author = {Lovelock, D. and Rund, H.},
    year = {1989},
    publisher = {Dover}
}

@article{lorenz1963,
  author = {Lorenz, Edward N.},
  title = {Deterministic Nonperiodic Flow},
  journal = {Journal of the Atmospheric Sciences},
  volume = {20},
  number = {2},
  pages = {130--141},
  year = {1963},
  doi = {10.1175/1520-0469(1963)020<0130:DNF>2.0.CO;2},
  publisher = {American Meteorological Society}
}

@book{Lyapunov1992,
  title = {The General Problem of the Stability of Motion},
  author = {Lyapunov, Aleksandr Mikhailovich},
  translator = {Fuller, A. T.},
  year = {1992},
  publisher = {Taylor and Francis},
  note = {English translation of the 1892 Russian original},
  isbn = {978-0750300744}
}

@article{Minkowski1908,
  author = {Minkowski, H.},
  title = {{Die Grundgleichungen für die elektromagnetischen Vorgänge in bewegten Körpern}},
  journal = {Nachrichten von der Gesellschaft der Wissenschaften zu Göttingen},
  volume = {1908},
  pages = {53-111},
  year = {1908},
  url = {http://eudml.org/doc/58707}
}

@article{RicciCurbastroLeviCivita1901,
  author = {Gregorio Ricci-Curbastro and Tullio Levi-Civita},
  title = {Méthodes de calcul différentiel absolu et leurs applications},
  journal = {Mathematische Annalen},
  volume = {54},
  pages = {125--201},
  year = {1901},
  doi = {10.1007/BF01456831},
  url = {https://doi.org/10.1007/BF01456831},
}

@inproceedings{Riemann1867,
  author = {Bernhard Riemann},
  title = {{Über die Hypothesen, welche der Geometrie zu Grunde liegen}},
  booktitle = {Gesammelte Werke},
  publisher = {Teubner},
  year = {1867},
  pages = {254--287},
  url = {https://de.wikipedia.org/wiki/Riemannscher_Krümmungstensor}
}

@article{Schuler1923,
  author = {Schuler, Maximilian},
  title = {{Die St\"orrung von Pendel- und Kreiselapparaten durch die Bewegung des Fahrzeuges}},
  journal = {Physikalische Zeitschrift},
  volume = {24},
  pages = {594--596},
  year = {1923}
}

@inproceedings{Schmalz2007Analytic,
  author = {Christian Schmalz and Winfried Lohmiller and Thomas Koehler},
  title = {Analytic Error Computation of the Strap-Down-Algorithm},
  booktitle = {AIAA Guidance, Navigation and Control Conference and Exhibit},
  year = {2007},
}

@book{slotine1991applied,
  title={Applied Nonlinear Control},
  author={Slotine, Jean-Jacques and Li, Weiping},
  year={1991},
  publisher={Prentice Hall},
  address={Englewood Cliffs, NJ}
}

@article{tabareau2007geometry,
  title={Geometry of the superior colliculus mapping and efficient oculomotor computation},
  author={Tabareau, Nicolas and Girard, Benoît and Bennequin, Daniel and Berthoz, Alain and Slotine, Jean-Jacques},
  journal={Biological Cybernetics},
  volume={97},
  number={4},
  pages={279--292},
  year={2007},
  doi={10.1007/s00422-007-0172-2}
}

@article{nguyen2017contraction,
  title={Contraction Analysis of Nonlinear DAE Systems},
  author={Nguyen, Hung D. and Slotine, Jean-Jacques},
  journal={arXiv preprint arXiv:1702.07421},
  year={2017},
  url={https://arxiv.org/abs/1702.07421}
}

@article{sontag,
  author  = {G. Russo and M. di Bernardo and E. D. Sontag},
  title   = {A contraction approach to the hierarchical analysis and design of networked systems},
  journal = {IEEE Transactions on Automatic Control},
  volume  = {58},
  number  = {5},
  pages   = {1328--1331},
  year    = {2013}
}

@book{titterton2004strapdown,
  title     = {Strapdown Inertial Navigation Technology},
  author    = {David H. Titterton and John L. Weston},
  year      = {2004},
  edition   = {2nd},
  publisher = {The Institution of Engineering and Technology},
  series    = {Radar, Sonar and Navigation},
  isbn      = {0863413587},
  address   = {London},
  pages     = {574}
}

@article{PartialContraction,
  author  = {Wang, Wei and Slotine, Jean-Jacques},
  title   = {On partial contraction analysis for coupled nonlinear oscillators},
  journal = {Biological Cybernetics},
  volume  = {92},
  number  = {1},
  pages   = {38--53},
  year    = {2005}
}

@article{BeyondConvexity,
  author  = {Patrick M. Wensing and Jean-Jacques Slotine},
  title   = {Beyond convexity---Contraction and global convergence of gradient descent},
  journal = {PLOS ONE},
  year    = {2020},
  volume  = {15},
  number  = {8},
  pages   = {e0236661},
  doi     = {10.1371/journal.pone.0236661}
}

\end{document}